\documentclass[11pt]{article}
\usepackage{theorem}
\usepackage{amssymb}
\usepackage{amsmath}
\usepackage{amsfonts}

\def \RR {\mathbb R}
\def \EE {\mathbb E}
\def \eps {\varepsilon}

\def \S {\mathcal S}

\newtheorem{theorem}{Theorem}[section]
\newtheorem{lemma}[theorem]{Lemma}

\newtheorem{proposition}[theorem]{Proposition}
\newtheorem{corollary}[theorem]{Corollary}
 {\theorembodyfont{\rmfamily}}

\begin{document}
\title{Rate of convergence of geometric symmetrizations}
\author{B. Klartag\thanks{Supported in part by the Israel
 Science Foundation and by the Minkowski center for Geometry.}
 \\ \\
School of Mathematical Sciences, \\
Tel Aviv University, \\
Tel Aviv 69978, Israel }
\date{}
\maketitle

\abstract{It is a classical fact, that given an arbitrary convex
body $K \subset \RR^n$, there exists an appropriate sequence of
Minkowski symmetrizations (or Steiner symmetrizations), that
converges in Hausdorff metric to a Euclidean ball. Here we provide
quantitative estimates regarding this convergence, for both
Minkowski and Steiner symmetrizations. Our estimates are
polynomial in the dimension and in the logarithm of the desired
distance to a Euclidean ball, improving previously known
exponential estimates. Inspired by a method of Diaconis \cite{D},
our technique involves spherical harmonics. We also make use of an
earlier result by the author regarding ``isomorphic Minkowski
symmetrization''. }

\section{Introduction}

Let $K \subset \RR^n$ be a convex body, and denote by $| \cdot |$
and $\langle \cdot, \cdot \rangle$ the usual Euclidean norm and
scalar product in $\RR^n$. Given a vector $u \in S^{n-1} = \{x \in
\RR^n ; |x| = 1\}$, we denote by $\pi_u(x) = x - 2 \langle x, u
\rangle u$ the reflection operator with respect to the hyperplane
through the origin, which is orthogonal to $u$ in $\RR^n$. The
result of a Minkowski symmetrization (sometimes called Blaschke
symmetrization) of $K$ with respect to $u$, is the body
$$ \tau_u(K) = \frac{K + \pi_u(K)}{2} $$
where the Minkowski sum of two sets $A,B \subset \RR^n$ is defined
as $A + B = \{ a + b ; a \in A, b \in B \}$. Let $h_K$ denote the supporting
functional of $K$, i.e. for $u \in \RR^n$
$$ h_K(u) = \sup_{x \in K} \langle x, u \rangle. $$
Then $h_{\tau_u(K)}(v) = \frac{1}{2} \left[ h_K(v) + h_K(\pi_u(v))
\right]$. The mean width of $K$ is defined as $w(K) = 2M^*(K) =
2\int_{S^{n-1}} h_K(u) d\sigma(u)$, where $\sigma$ is the unique
rotation invariant probability measure on the sphere. The mean
width is preserved under Minkowski symmetrizations.

Steiner symmetrization of $K$ with respect to a hyperplane $H$
yields the unique body $S_H(K)$ such that for any line $l$
perpendicular to $H$,
\renewcommand{\labelenumi}{(\roman{enumi})}
\begin{enumerate}
\item $S_H(K) \cap l$ is a closed segment whose center lies on $H$. \item
$Meas(K \cap l) = Meas (S_H(K) \cap l)$.
\end{enumerate}
\renewcommand{\labelenumi}{(\arabic{enumi})}
where $Meas$ is the one dimensional Lebesgue measure in the line
$l$. Steiner symmetrization preserves the volume of a set
and transforms convex sets to convex sets. See e.g. \cite{BF} for
more information about these symmetrizations, and their
applications in proving geometric inequalities.

\medskip
Consecutive Minkowski/Steiner symmetrizations may cause a convex
body to resemble a Euclidean ball. Starting with an arbitrary
convex body, one may apply a suitable sequence of
Minkowski/Steiner symmetrizations, and obtain a sequence of bodies
that converges to a Euclidean ball. This Euclidean ball would have
the same mean width/volume as had the original body.
 In this note, we investigate the rate
of this convergence. We ask how many symmetrizations are needed,
in order to transform an arbitrary convex body $K \subset \RR^n$
into a body that is $\eps$-close to a Euclidean ball. Our question
is ``almost isometric'' in its nature, as we try to provide
reasonable estimates even for small values of $\eps$. Previous
results in the literature are mostly of ``isomorphic'' nature, in
the sense that the symmetrization process is aimed at obtaining a
body which is uniformly ``isomorphic'' to a Euclidean ball (a body
is ``isomorphic'' to a Euclidean ball if its distance to a
Euclidean ball is bounded by some fixed, universal constant).

\smallskip The first quantitative result regarding Minkowski symmetrization
appears in \cite{BLM1}. Denote by $D$ the standard Euclidean ball
in $\RR^n$. Their result reads as follows:
\begin{theorem}
Let $0 < \eps < 1$, $n > n_0(\eps)$. Given an arbitrary convex
body $K \subset \RR^n$, there exist $c n\log n + c(\eps) n$
Minkowski symmetrizations that transform $K$ into a body $\tilde{K}$ such
that
$$ (1 - \eps) M^*(K) D \subset \tilde{K} \subset (1 + \eps) M^*(K)
D $$ where $c(\eps), n_0(\eps)$ are of the order of $exp(c \eps^{-2}
|\log \eps|)$ and $c > 0$ is a numerical constant.
\end{theorem}

Their proof uses the method of random Minkowski symmetrizations.
In \cite{K2}, the notion of randomness was altered, and has lead
to an improvement of the dependence on the dimension $n$. The
following is proved in \cite{K2}:
\begin{theorem}
Let $n \geq 2$ and let $K \subset \RR^n$ be a convex body. Then
there exist $5n$ Minkowski symmetrizations, such that when applied
to $K$, the resulting body $\tilde{K}$ satisfies,
$$ \left( 1 - c \frac{|\log \log n|}{\sqrt{\log n}} \right) M^*(K) D
\subset \tilde{K} \subset
\left( 1 + c \frac{|\log \log n| }{\sqrt{\log n}} \right) M^*(K) D
$$ where $c > 0$ is some numerical constant.
 \label{5n_theorem}
\end{theorem}
Note that both in \cite{K2} and in \cite{BLM1}, for any fixed
dimension, one cannot even formally conclude that there is
convergence to a Euclidean ball. This note fills that gap in the
literature, and also provides surprisingly good dependence on
$\eps$. The following theorem is proved here:
\begin{theorem}
Let $n \geq 2$, $0 < \eps < \frac{1}{2}$, and let $K \subset
\RR^n$ be a convex body. Then there exist $c n \log
\frac{1}{\eps}$ Minkowski symmetrizations, that transform $K$ into
a body $\tilde{K}$ that satisfies
$$ (1 - \eps) M^*(K) D \subset \tilde{K} \subset (1 + \eps) M^*(K)
D $$ where $c > 0$ is some numerical constant. \label{eps_theorem}
\end{theorem}

\medskip Our approach to the problem of Minkowski symmetrization
involves a number of novel ideas. First, rather than applying random Minkowski
symmetrizations, at each step we apply
$n$ symmetrizations with respect to the vectors of some random
orthonormal basis. This change of randomness improves the rate of
convergence by a factor of $\log n$ (see \cite{K1}, \cite{K2} and
also the remark following Corollary \ref{l_2_decay_not_harmonic}
here). Second, the use of spherical harmonics allows us to obtain
good estimates regarding symmetrization of polynomials on the
sphere. Finally, we approximate the supporting functional of $K$
with an appropriate polynomial (applying Theorem \ref{5n_theorem}
and a Jackson type theorem), and use the estimates obtained for
symmetrization of polynomials.

\bigskip Quantitative estimates regarding Steiner symmetrization
are more difficult to obtain, as the problem is non-linear. The
earliest estimate in the literature is due to Hadwiger \cite{H}.
It gives an estimate of the order of $\left(c
\frac{\sqrt{n}}{\eps^2} \right)^n$ for the number of Steiner
symmetrizations required in order to transform an arbitrary
$n$-dimensional convex body, to become $\eps$-close to a
Euclidean ball. In addition, an isomorphic result
appears in \cite{BLM2}, which was improved by a logarithmic factor
in \cite{KM}. The following is proved in \cite{KM}:
\begin{theorem}
Let $n \geq 2$ and let $K \subset \RR^n$ be a convex body, with
$Vol(K) = Vol(D)$. Then there exist $3n$ Steiner symmetrizations,
such that when applied to $K$, the resulting body $\tilde{K}$
satisfies, $$ cD \subset \tilde{K} \subset C D $$ where $c, C > 0$
are some numerical constants. \label{steiner_isomorphic}
\end{theorem}
Some related estimates also appear in \cite{T}. Our result is the
first estimate which is polynomial in $n$ and in $\log
\frac{1}{\eps}$. This shows that the precise geometric shape of a
convex body cannot prevent fast symmetrization of the body into an
almost Euclidean ball. In this note we shall prove the following
theorem.
\begin{theorem}
Let $K \subset \RR^n$ be a convex body, and let $0 < \eps <
\frac{1}{2}$. Let $r > 0$ be such that $Vol(K) = Vol(r D)$. Then there exist $c n^4 \log^2 \frac{1}{\eps}$
Steiner symmetrizations, that transform $K$ into a body
$\tilde{K}$ that satisfies
$$ (1 - \eps) r D \subset \tilde{K} \subset (1 + \eps) r
D $$ where $c > 0$ is some numerical constant. \label{steiner_eps}
\end{theorem}
The powers of $n$ and $\log \frac{1}{\eps}$ in Theorem
\ref{steiner_eps} seem non optimal. We conjecture that $c n \log
\frac{1}{\eps}$ Steiner symmetrizations are sufficient. Regarding
Minkowski symmetrizations, our result is tight in the sense
that the powers in Theorem \ref{eps_theorem} cannot be improved.

The proof of Theorem \ref{steiner_eps} is an application of
Theorem \ref{eps_theorem} and of a geometric result by Bokowski
and Heil. Throughout this paper, we denote by $c, C, c^{\prime}$
etc. positive numerical constants whose value is not necessarily
equal in different appearances.

\section{Spherical Harmonics}
\label{sphericals}

In this section we summarize a few facts about spherical
harmonics, to be used later on. For a comprehensive discussion on
the subject, we refer the reader to the concise expositions in
\cite{SW}, chapter $IV.2$, in \cite{M} and in \cite{G}. $P_k:
\RR^n \rightarrow \RR$ is a homegeneous harmonic of degree $k$, if
$P_k$ is a homogeneous polynomial of degree $k$ in $\RR^n$, and
$P_k$ is harmonic (i.e. $\triangle P_k \equiv 0$). We denote,
$$ \S_k = \{ P|_{S^{n-1}} \ ; \ P:\RR^n \rightarrow \RR \ is \ a \
homogenous \ harmonic \ of \ degree \ k \} $$ where $P|_{S^{n-1}}$
is the restriction of the polynomial $P$ to the sphere. $\S_k$ is the
space of spherical harmonics of degree $k$. It is
a linear space of dimension $\frac{(2k+n-2) \ (n+k-3)!}{k! \
(n-2)!}$. For $k \neq k^{\prime}$, the spaces $\S_k$ and
$\S_{k^{\prime}}$ are orthogonal to each other in $L_2(S^{n-1})$.
In addition, if $P$ is a polynomial of degree $k$ in $\RR^n$, then
$P|_{S^{n-1}}$ can be expressed as a sum of spherical harmonics of
degrees not larger than $k$. Therefore, $L_2(S^{n-1}) =
\bigoplus_k \S_k$. Spherical harmonics possess many symmetry
properties, partly due to their connection with the
representations of $O(n)$ (e.g. \cite{V}, chapter 9). For a fixed
dimension $n$, the Gegenbauer polynomials
$\{G_i(t)\}_{i=0}^{\infty}$ are defined by the following three
conditions:
\renewcommand{\labelenumi}{(\roman{enumi})}
\begin{enumerate}
\item $G_i(t)$ is a polynomial of degree $i$ in one variable.
\item For any $i \neq j$ we have $\int_{-1}^1 G_i(t) G_j(t) \left(
1 - t^2 \right)^{\frac{n-3}{2}} dt = 0 $. \item $G_i(1) = 1$ for
any $i$.
\end{enumerate}
\renewcommand{\labelenumi}{\arabic{enumi}.}
The Gegenbauer polynomials are closely related to spherical
harmonics. Next, we reformulate Lemma 3.5.4 from \cite{G}, which
is credited to Schneider. This useful lemma also follows from
Corollary 2.13, chapter $IV$ of \cite{SW}, and is true for all $n
\geq 2$.

\begin{lemma} Let $g \in \S_k$ be such that $\| g \|_2^2 = \int_{S^{n-1}} g^2(x) d\sigma(x) = 1$. Then,
$$ \int_{O(n)} g(U^{-1} x) g(U^{-1} y) d \mu(U) = G_k(\langle x, y \rangle) $$
where $\mu$ is the Haar probability measure on $O(n)$.
\label{legendre_lemma_corollary}
\end{lemma}

The following lemma reflects the fact that $\S_k$ is an
irreducible representation space of $O(n)$. We denote by $Proj_{\S_k}:L_2(S^{n-1}) \rightarrow \S_k$ the orthogonal projection onto $\S_k$.

\begin{lemma}
Let $f \in L_2(S^{n-1})$, and let $g \in \S_k$ be such that $\|g
\|_2 = 1$. Then,
\begin{equation}
\int_{O(n)} \left( \int_{S^{n-1}} f(Ux) g(x) d\sigma(x) \right)^2
d\mu(U) = \frac{\| Proj_{\S_k}(f) \|_2^2 }{dim(\S_k)}
\label{prod_harmonic}
\end{equation}
where $\mu$ is the Haar probability measure on $O(n)$.
\label{prod_harmonic_lemma}
\end{lemma}

\emph{Proof:}  Let $\{g_1,..,g_N\}$ be an orthonormal basis of
$\S_k$. Then,
\begin{equation}
\sum_{i=1}^{dim(\S_k)} \int_{O(n)} \left( \int_{S^{n-1}} f(Ux)
g_i(x) d\sigma(x) \right)^2 d\mu(U) \label{sum_basis_F_k}
\end{equation}
$$ = \int_{O(n)} \|
Proj_{\S_k}(f \circ U) \|_2^2 d\mu(U) = \| Proj_{\S_k}(f) \|_2^2
$$
because of the rotation invariance of $\S_k$. Therefore, it is
sufficient to prove that the integral in (\ref{prod_harmonic})
does not depend on the choice of $g \in \S_k$, as long as it
satisfies $\| g \|_2 = 1$. Indeed, in that case each of the
summands in (\ref{sum_basis_F_k}) equals $\frac{\| Proj_{\S_k}(f)
\|_2^2}{dim(\S_k)}$, for an arbitrary orthonormal basis
$\{g_1,..,g_N\}$ of $\S_k$. Let us try to simplify the integral in
(\ref{prod_harmonic}):
$$ \int_{O(n)} \int_{S^{n-1}} f(Ux) g(x) d\sigma(x) \int_{S^{n-1}}
f(Uy) g(y) d \sigma(y) d\mu(U) $$ $$ = \int_{S^{n-1}} \int_{S^{n-1}} f(x)
f(y) \int_{O(n)} g(U^{-1}x) g(U^{-1}y) d \mu(U) d\sigma(x)
d\sigma(y). $$
By Lemma \ref{legendre_lemma_corollary},
$\int_{O(n)} g(U^{-1}x) g(U^{-1}y) d \mu(U) = G_k(\langle x, y \rangle)$. Hence, the integral in
(\ref{prod_harmonic}) equals
$$ \int_{S^{n-1}} \int_{S^{n-1}} f(x) f(y) G_k(\langle x, y
\rangle) d\sigma(x) d\sigma(y) $$ which does not depend on $g$,
and the lemma is proved. \hfill $\square$

\section{Spherical Harmonics and Minkowski Symmetrization}

In this section we apply a series of Minkowski symmetrizations to
a convex body $K \subset \RR^n$. Each step in the symmetrization
process consists of symmetrizing $K$ with respect to the $n$
vectors of an orthonormal basis $\{e_1,..,e_n\}$ in $\RR^n$. Such
a step is denoted here as an ``orthogonal symmetrization'' with
respect to $\{e_1,..,e_n\}$. Applying an ``orthogonal
symmetrization'' with respect to $\{e_1,..,e_n\}$ to $K$, yields
a body denoted by $K^{\prime}$. Let $h$ be the supporting
functional of $K$, and $h^{\prime}$ be the supporting functional
of $K^{\prime}$. Then,
\begin{equation}
h^{\prime}(x) = \EE_{\eps} \ h \left(\sum_{i=1}^n \eps_i \langle x, e_i
\rangle e_i \right) \label{ortho_symm}
\end{equation}
 where the expectation is over $\eps \in \{\pm
1\}^n$, with respect to the uniform probability measure on the
discrete cube. Note that by (\ref{ortho_symm}), orthogonal
symmetrization may be viewed as an operation on support functions,
rather than on convex bodies. Furthermore, we may apply an
``orthogonal symmetrization'' to any function on the sphere, which
is not necessarily a support function of a convex body. Next, we
analyze the effect of orthogonal symmetrizations on spherical
harmonics.

\medskip Let $k$ be a positive integer. A function $g \in L_2(S^{n-1})$
is called ``invariant with respect to the orthonormal basis
$\{e_1,..,e_n\}$'', if for any $\eps \in \{\pm 1\}^n$, we have
$g(x) = g \left(\sum_i \eps_i \langle x, e_i \rangle e_i \right)$.
For a fixed orthonormal basis $\{e_1,..,e_n\}$ in $\RR^n$, we
denote by $\S_k^0$ the linear space of all invariant functions in
$\S_k$. Let $Proj_{\S_k^0}:\S_k \rightarrow \S_k^0$ be the
orthogonal projection in $L_2(S^{n-1})$. Then for $g \in \S_k$,
$$ g^{\prime}(x) = \EE_{\eps} \ g \left(\sum_{i=1}^n \eps_i \langle x, e_i
\rangle e_i \right) \ \ \ \Longleftrightarrow \ \ \ g^{\prime} =
Proj_{\S_k^0}(g), $$ i.e. the orthogonal symmetrization of $g$ is
the projection of $g$ onto $\S_k^0$.

\begin{lemma} If $k$ is odd, $dim(\S_k^0) = 0$. Otherwise,
$$ dim(\S_k^0) = \left( \! \! \! \begin{array}{c} n + \frac{k}{2} -
2 \\ n - 2 \end{array} \! \! \! \right). $$ \label{dim_F_K_0}
\end{lemma}

\vspace{-11pt} \emph{Proof:} The odd case is easy, since for $g
\in \S_k$ we necessarily have $g(x) = -g(-x)$, and for $g \in
\S_K^0$ we have $g(x) = g(-x)$. Hence, only $0 \in \S_k^0$. Next,
assume that $k$ is even, and let $g \in \S_k^0$ be an invariant
polynomial with respect to the basis $\{e_1,..,e_n\}$. We use the
coordinates $x_1,..,x_n$ with respect to this basis. Fixing
$x_2,..,x_n$ the polynomial $g$ satisfies $g_{x_2,..,x_n}(x_1) =
g_{x_2,..,x_n}(-x_1)$, and hence only even degrees of $x_1$ occur
in $g_{x_2,..,x_n}$. By repeating the argument for the rest of the
variables, we get that $g$ is a function of $x_1^2,..,x_n^2$
alone. We can write,
\begin{equation}
g(x_1,..,x_n) = \sum_{j=0}^{k/2} x_n^{2j} A_j(x_1,..,x_{n-1})
\label{one_var_out}
\end{equation}
where $A_j$ is a homogeneous polynomial of degree $k-2j$, which
depends solely on $x_1^2,..,x_{n-1}^2$. Let us calculate the
Laplacian of (\ref{one_var_out}):
$$ 0 = \sum_{j=1}^{k/2} 2j (2j -1 ) x_n^{2j-2}
A_j(x_1,..,x_{n-1}) + \sum_{j=0}^{k/2-1} x_n^{2j} \ \triangle
A_j(x_1,..,x_{n-1}) $$ or equivalently, $g \in \S_k^0$ if and anly
if for all $0 \leq j \leq \frac{k}{2}-1$,
\begin{equation}
(2j +2) (2j+1) A_{j+1} = -\triangle A_j. \label{coeff_cond}
\end{equation}
Therefore we are free to choose $A_0$ any way we like, as long as it
is a homogeneous polynomial of degree $k$, which involves only
even powers of the $n-1$ variables. When $A_0$ is fixed, $A_1,
A_2$ etc. are determined by equation (\ref{coeff_cond}), and
the function $g$ is recovered.

Hence, $dim(\S_k^0)$ equals the dimension of the space of the
possible $A_0(x_1^2,..,x_{n-1}^2)$, which is the dimension of the
space of all homogeneous polynomials of degree $k/2$ in $n-1$
variables. This number is known to be $ \left( \! \! \!
\begin{array}{c} n + \frac{k}{2} - 2 \\ n - 2 \end{array} \! \! \!
\right)$. \hfill $\square$

\medskip
We denote $N_k = dim(\S_k) = \left(\! \! \! \begin{array}{c} n+k-2
\\ n-2
\end{array} \! \! \! \right) \frac{n+2k-2}{n+k-2}$, and for an even $k$ denote
$N_k^0 = dim(\S_k^0) =\left(\! \! \! \begin{array}{c} n+k/2-2 \\
n-2
\end{array} \! \! \! \right)$. Clearly, these two quantities depend on $n$
which is absent from the notation, yet the
appropriate value of $n$ will be obvious from the context. We are
now ready to calculate the $L_2$ norm of a ``random orthogonal
symmetrization'' of a spherical harmonic - an orthogonal
symmetrization with respect to a basis that is chosen uniformly
over $O(n)$. Clearly, any ``orthogonal symmetrization'' of
an odd degree spherical harmonic vanishes. The even case is
treated in the following proposition.
\begin{proposition}
Let $k$ be a positive even integer, and let $g \in \S_k$ be a spherical
harmonic. We randomly select an orthonormal basis $\{ v_1,..,v_n
\} \in O(n)$, and symmetrize $g$ with respect to this basis. Then,
$$ \EE \| g^{\prime}_{v_1,..,v_n} \|_2^2 =  \frac{N_k^0}{N_k
} \ \| g \|_2^2 < \left( \frac{k}{n-2+k} \right)^{k/2} \| g
\|_2^2
$$ where the expectation is over the random choice of $\{
v_1,..,v_n \} \in O(n) $ (with respect to the Haar probability measure on $O(n)$).
\label{ortho_symm_l_2_decay}
\end{proposition}

\emph{Proof:} Fix an orthonormal basis $\{e_1,..,e_n\}$ of
$\RR^n$, and consider $\S_k^0$ with respect to that basis. Fix
also an orthonormal basis $S_1,..,S_{N_k^0}$ of $\S_k^0$. From the
discussion before Lemma \ref{dim_F_K_0},
$$ g^{\prime}_{e_1,..,e_n} = Proj_{\S_k^0}(g) $$
and if the columns of $U \in O(n)$ are $\{v_1,..,v_n\}$, then
$$  g^{\prime}_{v_1,..,v_n} = \left( Proj_{\S_k^0} (g \circ U)
\right) \circ U^{-1}. $$ Hence,
$$ \| g^{\prime}_{v_1,..,v_n} \|_2^2 = \| Proj_{\S_k^0} (g \circ U) \|_2^2 =
\sum_{j=1}^{N_k^0} \left(
\int_{S^{n-1}} g(Ux) S_j(x) d\sigma(x) \right)^2 $$ and by Lemma
\ref{prod_harmonic_lemma},
$$ \EE \| g^{\prime}_{v_1,..,v_n} \|_2^2 = \frac{\sum_{j=1}^{N_k^0}
\| g \|_2^2}{N_k} = \frac{N_k^0}{N_k} \| g \|_2^2. $$
 Note that
$$ \frac{N_k^0}{N_k} = \frac{n+k-2}{n+2k-2}
\prod_{i=1}^{k/2} \frac{(n+i-2)(k/2 + i)}{(n+2i-3)(n+2i-2)} <
\prod_{i=1}^{k/2} \frac{k/2 + i}{k/2+i + n - 2} $$
 which lies between
  $\left( \frac{k/2}{n-2+k/2} \right)^{k/2}$ and $\left( \frac{k}{n-2+k} \right)^{k/2}$.  \hfill $\square$

Since $\left( \frac{k}{n-2+k} \right)^{k/2}$ is a decreasing
function of $k$, then $\left( \frac{k}{n-2+k} \right)^{k/2} \leq
\frac{2}{n}$ for any $k \geq 2$, and we obtain the following
corollary:
\begin{corollary}
Let $f \in L_2(S^{n-1})$ satisfy $\int_{S^{n-1}} f(x) d\sigma(x) =
0$. We randomly select $\{ v_1,..,v_n \} \in O(n)$. Then,
$$ \EE \| f^{\prime}_{v_1,..,v_n} \|_2 < \frac{c}{\sqrt{n}} \|
f \|_2 $$ where the expectation is taken over the choice of $\{
v_1,..,v_n \} \in O(n)$, and $c = \sqrt{2}$.
\label{l_2_decay_not_harmonic}
\end{corollary}

\emph{Proof:} Expand $f$ into spherical harmonics: $f =
\sum_{k=1}^{\infty} f_k$ where $f_k = Proj_{\S_k}(f)$. Then
$f^{\prime}_{v_1,..,v_n} = \sum_{k=1}^{\infty} \left( f_k
\right)^{\prime}_{v_1,..,v_n}$ and
$$ \EE \| f^{\prime}_{v_1,..,v_n} \|_2^2 = \sum_{k=2}^{\infty}
\EE \| \left( f_k \right)^{\prime}_{v_1,..,v_n} \|_2^2 \leq \sum_k
\frac{2}{n} \| f_k \|_2^2 \leq \frac{2}{n} \|f\|_2^2. $$ An
application of Jensen inequality concludes the proof. \hfill
$\square$

\smallskip
\emph{Remark:} Using similar methods, one can prove that if $g \in
\S_k$ and $\tau_u(g)(x) = \frac{g(x) + g(\pi_u(x))}{2}$, then
$$ \EE_u \| \tau_u(g) \|_2^2 =  \frac{n-2+k}{n-2+2k} \| g \|_2^2. $$
Note the advantage of symmetrizing with respect to the $n$ vectors
of a random orthonormal basis, compared to symmetrization with
respect to $n$ random sphere vectors. For instance, if $k=2$ then
$$ \left( \frac{n-2+k}{n-2+2k} \right)^n \approx \frac{1}{e^2}. $$
Hence $n$ random symmetrizations may reduce the expectation of the
$L_2$ norm only by a constant factor.

\section{Decay of $L_{\infty}$ norm}

In Proposition \ref{ortho_symm_l_2_decay} and Corollary
\ref{l_2_decay_not_harmonic} we established a sharp estimate for
the decay of the $L_2$ norm under an ``orthogonal
symmetrization''. Now we deal with the more difficult problem of
estimating the decay of the $L_{\infty}$ norm of the function. Our
main tool is the following known lemma (see e.g. page 14 of
\cite{M}):
\begin{lemma}
Let $g \in \S_k$ be a spherical harmonic of degree $k$. Then,
$$\| g \|_{\infty} \leq \sqrt{ dim(\S_k) } \| g \|_2 =
\sqrt{N_k} \| g \|_2 $$ where $\| g \|_{\infty} = \sup_{x \in
S^{n-1}} |g(x)|$. \label{L_infinity}
\end{lemma}

We make use of the following well-known estimate of binomial
coefficients. For any $1 \leq k \leq n$,
\begin{equation}
\left( \frac{n}{k} \right)^k \leq \left(\! \! \! \begin{array}{c}
n \\ k
\end{array} \! \! \! \right) < \left( e \frac{n}{k} \right)^k.
\label{gluskin}
\end{equation}
In the following combinatorial lemmas, ``$\log$'' is to be
understood as the natural logarithm.
\begin{lemma}
Let $\eps > 0$, $n \geq 3$, and let $k \geq 2$ be an integer. Then,
$$ N_k^{c_1 \frac{1 + \log \left( 1 + \frac{2}{\eps} \right) }
{1 + \log \left( 1 + \frac{k}{n} \right) } } > \frac{n}{\eps^3} $$
where $c_1 > 0$ is some numerical constant. \label{triv_1}
\end{lemma}

\emph{Proof:} Denote $\alpha = k/n$.

\noindent \emph{Case 1:} $\alpha < 2$. In this case, $1 + \log
\left(1 + \frac{k}{n} \right) < 3$, and for $c_1 > 9$,
$$ N_k^{c_1 \frac{1 + \log \left( 1 + \frac{2}{\eps} \right) }
{1 + \log \left( 1 + \frac{k}{n} \right) } } > N_k^{3 + 3\log
\left( 1 + \frac{2}{\eps} \right)}
 > N_k \cdot N_k^{3 \log \left( 1 + \frac{2}{\eps} \right)}
 > \frac{N_k}{\eps^3} \geq \frac{n}{\eps^3}
$$
since for $k \geq 2$ we always have $N_k \geq n \geq 3$.

\noindent \emph{Case 2:} $\alpha \geq 2$. In this case, $1 + \log
\left( 1 + \frac{k}{n} \right) < 2 \log \left( 1 + \frac{k}{n-2}
\right) $. By (\ref{gluskin}), $N_k > \left( \frac{n+k-2}{n-2}
\right)^{n-2}$. For $c_1 > 6$,
$$ N_k^{c_1 \frac{1 + \log \left( 1 + \frac{2}{\eps} \right) }
{1 + \log \left( 1 + \frac{k}{n} \right) } } > \left( \left( 1 +
\frac{k}{n-2} \right)^{n-2} \right)^{\frac{3 + 3\log \left( 1 +
\frac{2}{\eps} \right) }{\log \left(1 + \frac{k}{n-2} \right) } }
$$ $$ = e^{3(n-2)} \left( 1 + \frac{2}{\eps} \right)^{3(n-2)} >
\frac{n}{\eps^3} $$ for any $n \geq 3$. \hfill $\square$

\begin{lemma}
Let $n \geq 3$, and let $k = \alpha n > 0$ be an even number. Then,
$$ \left( \frac{N_k^0}{N_k} \right)^T < \frac{1}{N_k} $$
for $T = c_2 \left[1 + \log (1+\alpha) \right] $, where $c_2 > 0$
is a numerical constant. \label{binomials}
\end{lemma}

\emph{Proof:} Since $\frac{n+2k-2}{n+k-2} > 1$, it is sufficient to
prove that
\begin{equation}
\left( \frac{\left( \! \! \! \begin{array}{c} n + k/2 - 2
\\ k/2 \end{array} \! \! \! \right)}{\left( \! \! \! \begin{array}{c} n + k - 2
\\ k \end{array} \! \! \!\right)} \right)^T < \frac{1}{\left( \! \! \! \begin{array}{c} n + k - 2
\\ k \end{array} \! \! \! \right)}.
\label{comb_1} \end{equation} \noindent \emph{Case 1:} $\alpha <
\frac{1}{2}$. The left hand side of (\ref{comb_1}) is equal to:
$$ \left( \prod_{i=1}^{k/2} \frac{k/2 + i}{k/2+i+n-2}
\right)^T < \left( \frac{k}{k+n-2} \right)^{\frac{kT}{2}} \leq
\left( \frac{k}{n} \right)^{\frac{kT}{2}}. $$ To obtain
(\ref{comb_1}) it is enough to prove that
$$ \left( \frac{k}{n} \right)^{\frac{kT}{2}} < \left( \frac{1}{e} \frac{k}{k+n-2} \right)^{k} $$
according to (\ref{gluskin}). Now, because $\alpha = \frac{k}{n} <\frac{1}{2}$, for $T=8$,
$$ \left( \frac{k}{n} \right)^{\frac{kT}{2}}
< \left( \frac{k}{n} \right)^{k} \left( \frac{1}{2}
\right)^{\frac{k(T-2)}{2}} <
 \left( \frac{2}{3e} \frac{k}{n} \right)^{k} < \left(
\frac{1}{e} \frac{k}{k+n-2} \right)^{k}. $$

\noindent \emph{Case 2:} $\alpha \geq \frac{1}{2}$. Since $\left(
\! \! \! \begin{array}{c} m \\ l \end{array} \! \! \! \right) =
\left( \! \! \! \begin{array}{c} m \\ m-l \end{array} \! \! \!
\right)$, the left hand side of (\ref{comb_1}) also equals:
$$ \left( \prod_{i=1}^{n-2} \frac{k/2 + i}{k+i}
\right)^T < \left( \frac{n-2+k/2}{n-2+k} \right)^{(n-2)T} < \left(
\frac{5}{6} \right)^{(n-2)T} $$ since $n-2 < 2k$ and because
$\frac{x+k/2}{x+k}$ is an increasing function of $x$. Now, for any
$T > \frac{1 + \log (1+\frac{k}{n-2})}{\log(6/5)}$,
$$ \left( \frac{5}{6} \right)^{(n-2)T}
< \left( \frac{1}{e} \frac{n-2}{n+k-2} \right)^{(n-2)} <
 \frac{1}{\left( \! \! \! \begin{array}{c} n + k - 2
\\ n-2 \end{array} \! \! \! \right)}. $$
Since for $n \geq 3$, we have $ \frac{1 + \log
(1+\frac{k}{n-2})}{\log(6/5)} < 10 \left[ 1 + \log ( 1 + \alpha)
\right]$, the lemma is proved.
 \hfill $\square$

\section{Proof of the Minkowski symmetrization result}

We make use of Jackson's theorem for the sphere, due to Newman and
Shapiro \cite{NS}:
\begin{theorem}
Let $n,k > 0$ be integers, and let $f:S^{n-1} \rightarrow \RR$ be
a $\lambda$-Lipschitz function on the sphere (i.e. $|f(x) - f(y)|
\leq \lambda |x-y|$ for any $x,y\in S^{n-1}$). Then there exists a
polynomial $P_k$ of degree $k$ in $n$ variables, such that for any
$x \in S^{n-1}$,
$$ | f(x) - P_k(x) | \leq c_3 \lambda \frac{n}{k} $$
\label{jackson} where $c_3 > 0$ is some numerical constant.
\end{theorem}

\emph{Proof of Theorem \ref{eps_theorem}:} We assume that $M^*(K)
= 1$. Begin with $5n$ symmetrizations, according to Theorem
\ref{5n_theorem}, to obtain a centrally-symmetric body $\bar{K}$.
Denote by $h$ its supporting functional. Then $h$ is a norm and
hence its Lipschitz constant equals $\sup_{x \in S^{n-1}} h(x)$. By
Theorem \ref{5n_theorem},
\begin{equation}
 \sup_{x \in S^{n-1}} h(x) < 1 + c \frac{|\log \log
n|}{\sqrt{\log n}} < c_4 \label{yap_5n_again}
\end{equation}
 for some numerical constant $c_4
> 0$. Hence $h$ is a $c_4$-Lipschitz function, and by Theorem
\ref{jackson}, there exists a polynomial $P_{\eps}(x)$ of degree $k
= \lceil \frac{n}{\eps} \rceil$ such that,
\begin{equation}
\sup_{x \in S^{n-1}} \left| P_{\eps}(x) - h( x) \right| < c_4 c_3
\eps. \label{jackson_five}
\end{equation}
Let $P_{\eps}(x) = \sum_{i=0}^k P_i(x)$ be the expansion of
$P_{\eps}$ into spherical harmonics. Randomly select $T$ orthonormal
bases (i.e. the bases are chosen independently and uniformly in
$O(n)$). Apply the corresponding $T$ orthogonal symmetrizations to
$P_{\eps}$ and $P_1,..,P_k$, to obtain the random polynomials
$P^{\prime}_{\eps}$ and $P^{\prime}_1,..,P^{\prime}_k$. Note that
still $P^{\prime}_{\eps} = \sum_{i=0}^k P^{\prime}_i$. Successive
application of Proposition \ref{ortho_symm_l_2_decay} yields that
for an even $i>0$,
$$ \EE \| P_i^{\prime} \|_2^2 = \left( \frac{N_i^0}
{N_i} \right)^T \ \| P_i \|_2^2. $$ Combining this with Lemma
\ref{L_infinity} (assume $n \geq 3$),
$$ \EE \| P_i^{\prime} \|_{\infty}^2 \leq N_i \left( \frac{N_i^0}
{N_i} \right)^T \ \| P_i \|_2^2. $$ Assume that $T > (c_1 + 1) c_2
\left[ 1 + \log \left( 1 + \frac{2}{\eps} \right) \right]$.
According to Lemma \ref{binomials},
$$ \EE \| P_i^{\prime} \|_{\infty}^2 < N_i \left( \frac{1}{N_i}
\right)^{(c_1 + 1) \frac{1 + \log \left(1+\frac{2}{\eps} \right)
}{1 + \log \left( 1 +  \frac{i}{n} \right) }} \ \| P_i \|_2^2
$$ $$ < N_i^{-c_1 \frac{1
+ \log \left(1+\frac{2}{\eps} \right)}{1 + \log \left( 1 +
\frac{i}{n} \right) } } \ \| P_i \|_2^2 < \frac{\eps^3}{n} \| P_i
\|_2^2 $$ where the last inequality follows from Lemma
\ref{triv_1}. Denote $I = P_0 = \int_{S^{n-1}} P_{\eps}(x)
d\sigma(x)$. Then,
$$ \EE \| P_{\eps}^{\prime}(x) - I \|_{\infty} \leq
\sum_{i=1}^{\left \lfloor \frac{k}{2} \right \rfloor} \EE \|
P_{2i}^{\prime}(x) \|_{\infty} \leq \sqrt{ \frac{k}{2}
\sum_{i=1}^{\left \lfloor \frac{k}{2} \right \rfloor} \EE \|
P_{2i}^{\prime} (x) \|_{\infty}^2 } $$
$$ \leq \sqrt{ \frac{1}{2} \left \lceil \frac{n}{\eps} \right
\rceil \sum_{i=1}^{\left \lfloor \frac{k}{2} \right \rfloor}
\frac{\eps^3}{n} \| P_{2i} (x) \|_2^2 } < \eps \| P_{\eps} \|_2 <
\eps \| P_{\eps} \|_{\infty} < \eps (c_4 + c_4 c_3 \eps)
$$
where the last inequality follows from (\ref{yap_5n_again}) and
(\ref{jackson_five}). Apply the same $T$ orthogonal
symmetrizations to $\bar{K}$, and obtain $K^{\prime}$. Denote by
$h^{\prime}$ the supporting functional of $K^{\prime}$. Then,
$$ \sup_{x \in S^{n-1}} \left| h^{\prime}(x) - P_{\eps}^{\prime}(x)
\right| < c_4 c_3 \eps, $$ and since by (\ref{jackson_five}) we
have $1- c_4 c_3 \eps < I < 1 + c_4 c_3 \eps$, then
$$ \EE \sup_{x \in S^{n-1}} \left| h^{\prime}(x) - 1 \right| <
c_4 c_3 \eps + c_4 c_3 \eps + \eps (c_4 + c_4 c_3 \eps) <
c^{\prime} \eps.
$$ Clearly,
$$ \sup_{x \in S^{n-1}} \left| h^{\prime}(x) - 1 \right| < c^{\prime} \eps \
\ \ \ \Rightarrow \ \ \ (1 - c^{\prime} \eps) D \subset K^{\prime}
\subset (1 + c^{\prime} \eps) D. $$ To summarize, we applied $5n +
(c_1 + 1) c_2 \left[ 1 + \log \left( 1 + \frac{2}{\eps} \right)
\right] n$ Minkowski symmetrizations to an arbitrary convex body,
some of which were chosen randomly. As a result of these
symmetrizations, we obtained a body such that the expectation of
its distance to a Euclidean ball is no more than $c^{\prime} \eps$.
Therefore, there exists some numerical constant $c > 0$, and $c n
\log \frac{1}{\eps}$ symmetrizations that bring the body to be
$\eps$-close to a Euclidean ball. \hfill $\square$

\emph{Remarks:}
\begin{enumerate}
\item The case $n=2$ should be treated separately. In this case,
$dim(\S_k) = 2, dim(\S_k^{\circ}) = 1$ for any $k$. It is easy to
verify that the proof works in this case as well. \item Theorem
\ref{eps_theorem} is optimal in the sense that one cannot obtain
an estimate for the number of minimal symmetrizations, of the form
$f(n) g(\eps)$ with $f(n) << n$ or $g(\eps) << \log
\frac{1}{\eps}$. Indeed, the dependence on $n$ should be at least
linear, as it takes a segment $n-1$ symmetrizations just to become
$n$-dimensional. Regarding the dependence on $\eps$, if we take a
segment and apply any $\lfloor c \log \frac{1}{\eps} \rfloor$
symmetrizations, then the segment is transformed into a zonotope
which is a sum of no more than $ \frac{1}{\eps^c}$ segments. Even
in dimension two, this zonotope cannot be $\eps$-close to a
Euclidean ball, for a small enough $c$. \item Note that Theorem
\ref{eps_theorem} is not tight for all possible values of $n$ and
$\eps$. For example, Theorem \ref{5n_theorem} is better than
Theorem \ref{eps_theorem} when $\eps = c \frac{|\log \log
n|}{\sqrt{\log n}}$.
\end{enumerate}

\section{Application to Steiner Symmetrization}
In this section we prove Theorem \ref{steiner_eps}. We make use of
a result due to Bokowski and Heil. The following theorem is a
special case of Theorem 2 in \cite{BH} (the case $(i,j,k) =
(0,d-1,d)$ in the notations of that paper).
\begin{theorem} Let $K \subset R D$ be a convex body. Then,
$$ n^2 R^{n-1} M^*(K) \leq \frac{Vol(K)}{Vol(D)} + (n^2 - 1) R^n. $$
\label{boko_heil}
\end{theorem}

\vspace{-11pt} An immediate corollary follows:
\begin{corollary} Let $\eps >0$, and let $K \subset (1 + \eps) D$ be a convex body
in $\RR^n$ with $Vol(K) = Vol(D)$. Then,
$$ M^*(K) < 1 + \left(1 - \frac{1}{n^2} \right) \eps. $$
In addition, if $\eps < \frac{1}{n}$ then,
$$ M^*(K) < 1 + \left(1 - \frac{1}{2n} \right) \eps. $$
\label{mean_width_decay}
\end{corollary}

\vspace{-11pt} \emph{Proof}: By Theorem \ref{boko_heil}, since
$\frac{Vol(K)}{Vol(D)} = 1$,
$$ M^*(K) \leq (1 + \eps) \left( 1 - \frac{1}{n^2} \right) + \frac{1}{n^2(1 + \eps)^{n-1}}
< (1 + \eps) \left( 1 - \frac{1}{n^2} \right) + \frac{1}{n^2} $$
and therefore $ M^*(K) < 1 + \left(1 - \frac{1}{n^2} \right) \eps
$. Now, assume that $\eps < \frac{1}{n}$. Using the elementary
inequality $\frac{1}{(1 + \eps)^{n-1}} < 1 - (n-1) \eps +
\frac{n(n-1)}{2} \eps^2$, we obtain
$$ M^*(K) \leq (1 + \eps) \left( 1 - \frac{1}{n^2} \right) +
\frac{1}{n^2} \left[ 1 - (n-1) \eps + \frac{n(n-1)}{2} \eps^2
\right] $$ $$ < 1 + \eps - \frac{\eps}{n} + \frac{\eps^2}{2} < 1 +
\eps - \frac{\eps}{2n}. $$  \hfill $\square$

Given a convex body $K \subset \RR^n$, define $ R(K) = \inf \{ R
> 0; K \subset R D \}$.

\begin{lemma}
Let $K \subset \RR^n$ be a convex body with $Vol(K) = Vol(D)$.
Assume that there exists $0 < \eps < C$ such that $R(K) = 1 +
\eps$, where $C > 1$. Then there exist $c_5 n (\log \frac{1}{\eps} + \log n)$
Steiner symmetrizations that transform $K$ into $\tilde{K}$ such
that
$$ R(\tilde{K}) < 1 + \left(1 - \frac{1}{2n^2} \right) \eps $$
and if $\eps < \frac{1}{n}$,
$$ R(\tilde{K}) < 1 + \left(1 - \frac{1}{4n} \right) \eps $$
where $c_5 = c_5(C) > 0$ depends solely on $C$.
\label{essential_steiner}
\end{lemma}

\emph{Proof:} Let $\tilde{K}$ be the body obtained from $K$ after
the $c n \log \frac{4Cn^3}{\eps}$ symmetrizations given by Theorem
\ref{eps_theorem}. Despite the fact that Theorem \ref{eps_theorem}
is concerned with Minkowski symmetrizations, we apply the
corresponding Steiner symmetrization (with respect to the same
hyperplanes). Since Steiner symmetrizations are contained in
Minkowski symmetrizations,
$$ R(\tilde{K}) < \left( 1 + \frac{\eps}{4Cn^3} \right) M^*(K). $$
Apply corollary \ref{mean_width_decay} and the fact that $\eps <
C$ to get that
$$ R(\tilde{K}) < \left( 1 + \frac{\eps}{4Cn^3} \right)
 \left[ 1 + \left(1 - \frac{1}{n^2} \right) \eps \right] < 1 + \left( 1 - \frac{1}{2 n^2} \right) \eps $$
and if $\eps < \frac{1}{n}$,
$$ R(\tilde{K}) < \left( 1 + \frac{\eps}{4Cn^3} \right)
 \left[ 1 + \left(1 - \frac{1}{2n} \right) \eps \right] < 1 + \left( 1 - \frac{1}{4n} \right) \eps. $$
\hfill $\square$

\begin{proposition}
Let $n \geq 2$, $0 < \eps < \frac{1}{2}$, and let $K \subset
\RR^n$ be a convex body with $Vol(K) = Vol(D)$. Then there exist
$c_6 \left[ n^3 \log^2 n + n^2 \log^2 \frac{1}{\eps} \right]$
Steiner symmetrizations, that transform $K$ into $\tilde{K}$ which
satisfies
$$ R(\tilde{K}) < 1 + \eps $$
where $c_6 > 0$ is a numerical constant. \label{steiner_above}
\end{proposition}

\emph{Proof:} First, apply $3n$ Steiner symmetrizations to $K$,
according to Theorem \ref{steiner_isomorphic}, to obtain an
isomorphic Euclidean ball $\bar{K}$. Then,
$$ \bar{K} \subset C D. $$
Let us define a sequence of convex bodies: $K_0 = \bar{K}$, and
$K_i$ is obtained from $K_{i-1}$ using $c_5 n (\log
\frac{1}{R(K_{i-1}) - 1} + \log n)$ Steiner symmetrizations, as in
Lemma \ref{essential_steiner}. Then,
\begin{equation}
R(K_i) - 1  < \left(1 - \frac{1}{2n^2} \right) \left[ R(K_{i-1}) -
1 \right] < \left(1 - \frac{1}{2n^2} \right)^i \left[ R(K_0) - 1
\right]. \label{geometric_progression}
\end{equation}
Let $T_1$ be the minimal integer such that $$ R(K_{T_1}) < 1 +
\frac{1}{n}. $$ Since $R(K_0) < C$, then by
(\ref{geometric_progression}) necessarily $T_1 < c n^2 \log n$.
For any $i \leq T_1$ we have $R(K_{i-1}) \geq 1 + \frac{1}{n}$ and
hence by Lemma \ref{essential_steiner} we used no more than
$c^{\prime} n \log n$ symmetrizations to obtain $K_i$ from
$K_{i-1}$. In total, we used less than $\tilde{c} n^3 \log^2 n$
symmetrizations to obtain $K_{T_1}$. By Lemma
\ref{essential_steiner} for any $i > 0$,
$$ R(K_{{T_1} + i}) - 1< \left(1 - \frac{1}{4n} \right) \left[
R(K_{{T_1}+i-1}) - 1 \right] < \left(1 - \frac{1}{4n} \right)^i.
$$
Let $T_2$ be the first integer such that
$$ R(K_{T_1 + T_2}) < 1 + \eps. $$
Then $T_2 < c n \log \frac{1}{\eps}$. Define $\tilde{K} =
K_{T_1+T_2}$. For any $T_1 < i \leq T_1+T_2$ we used no more than
$c^{\prime} n (\log \frac{1}{\eps} + \log n)$ symmetrizations to
obtain $K_i$ from $K_{i-1}$. In total we applied a maximum of
$\tilde{c} n^3 \log^2 n + \tilde{c} n^2 \log^2 \frac{1}{\eps}$
Steiner symmetrizations. \hfill $\square$

Proposition \ref{steiner_above} proves the existence of a rather
small circumscribing ball for the symmetrized body. In order to
symmetrize the body from below, we use the following standard
lemma. Its proof is outlined for completeness.
\begin{lemma}
Let $0 < \eps < 1$, and let $K \subset \RR^n$ be a convex body
with $M^*(K) \geq 1$. Assume that $K \subset \left[1 + (c_7
\eps)^n \right] D$. Then
$$ (1 - \eps) D \subset K $$
where $c_7 > 0$ is some numerical constant. \label{small_cap}
\end{lemma}

\emph{Proof:} Assume on the contrary that there exists $x_0 \in
S^{n-1}$ with $\| x_0 \|_* < 1 - \eps$, where $\| \cdot \|_* =
h_K(\cdot)$. Then for $x \in S^{n-1}$ with $|x - x_0| <
\frac{\eps}{4}$, we have
$$ \| x \|_* \leq \| x_0 \|_* + \| x - x_0 \|_* < 1 - \eps + (1 +
(c_7 \eps)^n) |x - x_0| < 1 - \frac{\eps}{2}. $$ Denote $A = \{ x
\in S^{n-1} ; |x - x_0| \leq \frac{\eps}{4} \}$. Then,
$$ M^*(K) = \int_{S^{n-1}} \| x \|_* d\sigma(x) < (1 - \sigma(A)) \left(1 +
(c_7 \eps)^n \right) + \sigma(A) \left( 1 - \frac{\eps}{2} \right).
$$ The projection of $A$ onto the hyperplane orthogonal to $x_0$
contains a Euclidean ball of radius larger than $\frac{\eps}{4
\sqrt{2}}$. Therefore,
$$ \sigma(A) > \frac{Vol(D_{n-1})}{Vol(S^{n-1})} \left(
\frac{\eps}{4 \sqrt{2}} \right)^{n-1} > \frac{1}{\sqrt{\pi} n} \left(\frac{\eps}{4 \sqrt{2}}
\right)^{n-1} > \left( \frac{\eps}{30} \right)^{n-1}
$$ where $D_{n-1}$ is the $n-1$ dimensional Euclidean unit ball,
and $Vol$ is interpreted here as the $n-1$ dimensional volume.
Thus,
$$ M^*(K) < 1 + (c_7\eps)^n - \sigma(A) \frac{\eps}{2} < 1 + (c_7\eps)^n
- \left(\frac{\eps}{30}\right)^n $$ and for $c_7 = \frac{1}{30}$ we
obtain a contradiction. \hfill $\square$

\emph{Proof of Theorem \ref{steiner_eps}:} It is sufficient to consider
the case $Vol(K) = Vol(D)$. Apply Proposition
\ref{steiner_above} with $\eps^{\prime} = (c_7 \eps)^n$. We use
$$ c_6 \left[ n^3 \log^2 n + n^2 \log^2 \frac{1}{\eps^{\prime}} \right] < c^{\prime} n^4 \log^2 \frac{1}{\eps}
$$ Steiner symmetrizations, and obtain a body $\tilde{K}$ such that
$$ \tilde{K} \subset \left( 1 + (c_7 \eps)^n \right) D. $$
Since $Vol(K) = Vol(D)$, by Urysohn $M^*(K) \geq 1$. Using Lemma
\ref{small_cap},
$$ (1 - \eps) D \subset \tilde{K} \subset \left( 1 + (c_7 \eps)^n \right)
D \subset (1 + \eps) D $$ and the theorem is proved. \hfill
$\square$

\smallskip
\emph{Remark:} Theorem \ref{5n_theorem} is crucial to the proof of
Theorem \ref{eps_theorem}. Only after obtaining the precise
isomorphic statement regarding Minkowski symmetrization, can we
prove the sharp almost isometric version. However, in the proof of
Theorem \ref{steiner_eps} we may apply weaker estimates than that
in Theorem \ref{steiner_isomorphic}, and derive the same
conclusion. This could be another indication that the powers in
Theorem \ref{steiner_eps} are not optimal.

\bigskip \emph{Acknowledgement.} Part of the research was done
during my visit to the University of Paris $VI$, and I am grateful
for their hospitality. Thanks also to the anonymous referee for the thorough review of the
paper.

\end{document}